\documentclass[11pt]{article}
\usepackage[utf8]{inputenc}
\usepackage{amsmath}
\usepackage{amsthm}
\usepackage{amssymb}
\usepackage{enumerate}
\usepackage{hyperref}
\usepackage{aliascnt}
\usepackage{graphicx}
\usepackage{color}

\theoremstyle{plain}
\newtheorem{thm}{Theorem}[section]
\newaliascnt{lem}{thm}
\newtheorem{lem}[lem]{Lemma}
\aliascntresetthe{lem}
\newaliascnt{pro}{thm}
\newtheorem{pro}[pro]{Proposition}
\aliascntresetthe{pro}
\newaliascnt{cor}{thm}
\newtheorem{cor}[cor]{Corollary}
\aliascntresetthe{cor}
\newaliascnt{que}{thm}
\newtheorem{que}[que]{Question}
\aliascntresetthe{que}
\newaliascnt{con}{thm}
\newtheorem{con}[con]{Conjecture}
\aliascntresetthe{con}
\newaliascnt{clm}{thm}

\aliascntresetthe{clm}

\theoremstyle{definition}
\newaliascnt{exm}{thm}

\aliascntresetthe{exm}
\theoremstyle{plain}

\newcommand{\N}{\mathbb{N}}

\newcommand{\B}{\mathcal{B}}
\newcommand{\C}{\mathcal{C}}
\newcommand{\I}{\mathcal{I}}
\renewcommand{\L}{\mathcal{L}}

\newcommand{\R}{\mathcal{R}}

\DeclareMathOperator{\cl}{cl}
\DeclareMathOperator{\fin}{fin}
\DeclareMathOperator{\rank}{rk}

\begin{document}
\title{\scshape On the intersection of infinite matroids}
\author{Elad Aigner-Horev\footnote{Research supported by the Minerva foundation.} \and Johannes Carmesin \and Jan-Oliver Fröhlich}
\date{University of Hamburg\\9 July 2012}
\maketitle

\begin{abstract}
\noindent We show that the {\sl infinite matroid intersection conjecture} of Nash-Williams implies the infinite Menger theorem proved recently by Aharoni and Berger. 

We prove that this conjecture is true whenever one matroid is nearly finitary and the second is the dual of a nearly finitary matroid, where the nearly finitary matroids form a superclass of the finitary matroids.

In particular, this proves the infinite matroid intersection conjecture for finite-cycle matroids of $2$-connected, locally finite graphs with only a finite number of vertex-disjoint rays.
\end{abstract}

\section{Introduction} 

The infinite Menger theorem\footnote{see \autoref{thm:menger} below.} was conjectured by Erd\H{o}s in the $1960$s and proved recently by Aharoni and Berger~\cite{AharoniBerger}. 
It states that for any two sets of vertices $S$ and $T$ in a connected graph,
there is a set of 
vertex-disjoint $S$-$T$-paths whose maximality is witnessed by an $S$-$T$-separator picking exactly one vertex form each of these paths.

The complexity of the only known proof of this theorem and the fact that the finite Menger theorem has a short matroidal proof, make it natural to ask whether a matroidal proof of the infinite Menger theorem exists. In this paper, we propose a way to approach this issue by proving that a certain conjecture of Nash-Williams regarding infinite matroids implies the infinite Menger theorem. 

Recently, Bruhn, Diestel, Kriesell, Pendavingh and Wollan~\cite{matroid_axioms} 
found axioms for infinite matroids in terms of independent sets, bases, circuits, closure and
(relative) rank. These axioms allow for duality of infinite matroids as known
from finite matroid theory, which settled an old problem of Rado. 
With these new axioms it is possible now to look which theorems of finite
matroid theory have infinite analogues.

Here, we shall look at the \emph{matroid intersection theorem}, which is a classical result in finite matroid theory~\cite{Oxley}. It asserts that {\sl the maximum size of a common independent set of two matroids $M_1$ and $M_2$ on a common ground set $E$ is given by}
\begin{equation}\label{eqn:intersection}
\min_{X \subseteq E} \rank_{M_1}(X) + \rank_{M_2}(E \setminus X),
\end{equation}
where $\rank_{M_i}$ denotes the rank function of the matroid $M_i$.

In this paper, we consider the following conjecture of Nash-Williams, which
first appeared in~\cite{Aharoni:Ziv:92} and serves as an infinite analogue
to the finite matroid intersection theorem\footnote{An alternative notion of
infinite matroid intersection was recently proposed by
Christian~\cite{christian_phd}.}. 

\begin{con}\label{thm:nash-williams-aharoni}\emph{[The infinite matroid intersection conjecture]}\\
Any two matroids $M_1$ and $M_2$ on a 
common ground set $E$ have a common independent set $I$ admitting a partition $I
= J_1 \cup J_2$ such that $\cl_{M_1}(J_1) \cup \cl_{M_2}(J_2) = E$.
\end{con}

\noindent
Here, $\cl_M(X)$ denotes the \emph{closure} of a set $X$ in a matroid $M$; it
consists of $X$ and the elements spanned by $X$ in $M$
(see~\cite{Oxley}).
Originally, Nash-Williams's Conjecture just concerned \emph{finitary
matroids}, those all of whose circuits are finite.

\subsection{Our results}

Aharoni and Ziv~\cite{Aharoni:Ziv:92} proved that \autoref{thm:nash-williams-aharoni} implies the infinite analogues of 
König's and Hall's theorems. We strengthen this by showing that this
conjecture implies the celebrated {\sl infinite Menger theorem}
(\autoref{thm:menger} below), which is known to imply the 
infinite analogues of  König's and Hall's theorems \cite{DiestelBook10}.

\begin{thm}\label{thm:inter-to-menger}
The infinite matroid intersection conjecture for finitary matroids implies the infinite Menger theorem.
\end{thm}

In finite matroid theory, an exceptionally short proof of the matroid intersection theorem employing the well-known \emph{finite matroid union theorem}~\cite{Oxley,schrijverBook} is known. The latter theorem asserts that for two finite matroids $M_1 = (E_1,\I_1)$ and $M_2=(E_2,\I_2)$
the set system 
\begin{equation}\label{eqn:what-is-union}
\I(M_1\vee M_2) = \{I_1\cup I_2\mid I_1 \in \I_1,\; I_2 \in \I_2\}
\end{equation}
forms the set of independent sets of their \emph{union matroid} $M_1 \vee M_2$.

In a previous paper~\cite[Proposition~1.1]{union}, we showed that for infinite
matroids $M_1$ and $M_2$, the set system $\I(M_1\vee M_2)$ is not necessarily a
matroid. This then raises the question of whether the traditional connection
between (infinite) matroid union and intersection still holds. 
In this paper, we prove the following.

\begin{thm}\label{thm:intersection}
If $M_1$ and $M_2$ are matroids on a common ground set $E$ and $M_1 \vee M^*_2$
is a matroid, then \autoref{thm:nash-williams-aharoni} holds for $M_1$ and $M_2$.
\end{thm} 

\noindent
Throughout, $M^*$ denotes the dual of a matroid $M$. 

In \cite{union} we show that the `largest' class of matroids for which one can
have a union theorem is essentially a certain superclass of the finitary
matroids called the \emph{nearly finitary} matroids (to be defined next). 
This, together with   
\autoref{thm:intersection}, enables us to make additional progress on
\autoref{thm:nash-williams-aharoni}, as set out below.  

Nearly finitary matroids are defined as follows \cite{union}. 
For any matroid $M$, taking as circuits only the finite circuits of $M$ defines a (finitary) matroid with the same ground set as $M$. This matroid is called the \emph{finitarization} of $M$ and denoted by $M^{\fin}$.

It is not hard to show that every basis $B$ of $M$ extends to a
basis $B^{\fin}$ of $M^{\fin}$, 
and conversely every basis $B^{\fin}$ of $M^{\fin}$ contains a basis $B$ of $M$.
Whether or not $B^{\fin}\setminus B$ is finite will in general depend on the
choices for $B$ and $B^{\fin}$, but given a choice for one of the two, it will
no longer depend on the choice for the second one.

We call a matroid $M$ \emph{nearly finitary} if every base of its finitarization contains a base of $M$ such that their difference is finite. 

Next, let us look at some examples of nearly finitary matroids.
 There are three natural extensions to the notion of a finite graphic matroid in
an infinite context~\cite{matroid_axioms};
each with ground set $E(G)$. 
The most studied one is the \emph{finite-cycle matroid}, denoted $M_{FC}(G)$,  whose circuits are the finite cycles of $G$.
This is a finitary matroid, and hence is also nearly finitary.

The second extension is the \emph{algebraic-cycle matroid}, denoted $M_A(G)$, whose circuits are the finite cycles and double rays of 
$G$~\cite{{matroid_axioms}, {RD:HB:graphmatroids}}\footnote{$M_A(G)$ is not necessarily a matroid for any $G$; see~\cite{Higgs:axioms}.}.

\begin{pro}\label{thm:alg}
$M_A(G)$ is a nearly finitary matroid if and only if $G$ has only a finite number of vertex-disjoint rays.
\end{pro}  

The third extension is the \emph{topological-cycle matroid}, denoted $M_C(G)$\footnote{$M_C(G)$ is a matroid for any $G$;
see~\cite{RD:HB:graphmatroids}.}, whose circuits are the topological cycles of
$G$ (Thus $M_C^{\fin}(G)=M_{FC}(G)$ for any finitely separable graph $G$; see
\autoref{sec:topo} or~\cite{RD:HB:graphmatroids} for definitions).

\begin{pro}\label{thm:topo}
Suppose that $G$ is $2$-connected and locally finite.
Then, $M_C(G)$ is a nearly finitary matroid if and only if $G$ has only a finite number of vertex-disjoint rays.
\end{pro}

Having introduced nearly finitary matroids, we now state the result of~\cite{union}.

\begin{thm}\label{thm:nearly-union}\emph{[Nearly finitary union theorem~\cite{union}]}\\
If $M_1$ and $M_2$ are nearly finitary matroids, then $M_1 \vee M_2$ is a nearly finitary matroid.
\end{thm}

\noindent
The following is a consequence of \autoref{thm:nearly-union} and
\autoref{thm:intersection}.  

\begin{cor}\label{thm:nearly-intersection}
\autoref{thm:nash-williams-aharoni} holds for $M_1$ and $M_2$ whenever $M_1$ is nearly finitary and $M_2$ is the dual of a nearly finitary matroid.  
\end{cor}

Aharoni and Ziv~\cite{Aharoni:Ziv:92} proved that the infinite matroid intersection conjecture is true whenever one matroid is finitary and the other is a countable direct sum of finite-rank matroids.
Note that \autoref{thm:nearly-intersection} does not imply this result of~\cite{Aharoni:Ziv:92} nor is it implied by it.

\autoref{thm:topo} and \autoref{thm:nearly-intersection} can be used to prove the following.
\vspace{\baselineskip}

\begin{cor}\label{thm:int_MFC}
Suppose that $G$ and $H$ are $2$-connected, locally finite graphs with only a finite number of vertex-disjoint rays.
Then their finite-cycle matroids $M_{FC}(G)$ and $M_{FC}(H)$ satisfy the
intersection conjecture.
\end{cor}

Similar results are true for \emph{the algebraic-cycle matroid}, \emph{the
topological-cycle matroid}, and their duals.

This paper is organized as follows. Additional notation, terminology, and basic lemmas are given in \autoref{sec:pre}. In \autoref{sec:menger} we prove \autoref{thm:inter-to-menger}. In \autoref{sec:union-to-inter} we prove \autoref{thm:intersection}, and in \autoref{sec:graphic} we prove Propositions~\ref{thm:alg} and~\ref{thm:topo} and \autoref{thm:int_MFC}.

\section{Preliminaries}\label{sec:pre}

Notation and terminology for graphs are that of~\cite{DiestelBook10}, and for matroids  that of~\cite{Oxley,matroid_axioms}.

Throughout, $G$ always denotes a graph where $V(G)$ and $E(G)$ denote its vertex and edge sets, respectively. 
We write $M$ to denote a matroid and write $E(M)$, $\I(M)$, $\B(M)$, and $\C(M)$ to denote its ground set, independent sets, bases, and circuits, respectively. 

We review the definition of a matroid as this is given in~\cite{matroid_axioms}.
A set system $\I$ is the set of independent sets of a matroid if it satisfies the following \emph{independence axioms}:
\begin{itemize}
	\item[(I1)] $\emptyset\in \I$.
	\item[(I2)] $\left\lceil \I \right\rceil=\I$, that is, $\I$ is closed under taking subsets.
	\item[(I3)] Whenever $I,I'\in \I$ with $I'$ maximal and $I$ not maximal, there exists an $x\in I'\setminus I$ such that $I+x\in \I$.
	\item[(IM)] Whenever $I\subseteq X\subseteq E$ and $I\in\I$, the set $\{I'\in\I\mid I\subseteq I'\subseteq X\}$ has a maximal element.
\end{itemize}

The following is a well-known fact for finite matroids (see, e.g., \cite{Oxley}), which generalizes easily 
to infinite matroids.  
\begin{lem}\label{thm:meet}\emph{\cite[Lemma 3.11]{matroid_axioms}}\\
Let $M$ be a matroid. Then, 
$
|C \cap C^*| \not= 1,
$
whenever $C \in \C(M)$ and $C^* \in \C(M^*)$.
\end{lem}

We end this section with the definition of {\sl exchange chains}. 
For a set $X \subseteq E(M)$, an $X$-\emph{circuit} is a circuit containing $X$. 
For sets $I_1\in\I(M_1)$ and $I_2\in\I(M_2)$, and elements $x\in I_1\cup I_2$ and  $y\in E(M_1)\cup E(M_2)$ (possibly in $I_1\cup I_2$), a tuple $Y=(y_0=y, \ldots, y_n=x)$ is called an \emph{even $(I_1, I_2,y,x)$-exchange chain} (or \emph{even $(I_1, I_2,y,x)$-chain}) of \emph{length} $n$ if the following terms are satisfied. 

\begin{enumerate}[(X1)]
	\item For an even $i$, there exists a $\{y_i, y_{i+1}\}$-circuit $C_i\subseteq I_1+y_i$ of $M_1$.
	\item For an odd $i$, there exists a $\{y_i, y_{i+1}\}$-circuit $C_i\subseteq I_2+y_i$ of $M_2$.
\end{enumerate}

\noindent
If $n\geq 1$, then (X1) and (X2) imply that $y_0\notin I_1$ and that, starting with $y_1\in I_1\setminus I_2$, the elements $y_i$ alternate between $I_1\setminus I_2$ and $I_2\setminus I_1$; the single exception being $y_n$ which might lie in $I_1\cap I_2$.

By an \emph{odd exchange chain} (or \emph{odd chain}) we mean an even chain with the words `even' and `odd' interchanged in the definition.
Consequently, we say \emph{exchange chain} (or \emph{chain}) to refer to either of these notions.
Furthermore, a subchain of a chain is also a chain;
that is, given an $(I_1, I_2, y_0, y_n)$-chain $(y_0,\ldots, y_n)$, the tuple $(y_k,\ldots, y_l)$ is an $(I_1, I_2, y_k, y_l)$-chain for $0\leq k\leq l\leq n$.

\begin{lem}\label{thm:chain}\emph{\cite[Lemma 4.4]{union}}\\
If there exists an $(I_1, I_2,y,x)$-chain, then $(I + y) - x\in\I(M_1\vee M_2)$ where $I := I_1\cup I_2$.
Moreover, if $x\in I_1\cap I_2$, then $I+y\in\I(M_1\vee M_2)$.
\end{lem}

\section{From infinite matroid intersection to the infinite Menger theorem}\label{sec:menger}

In this section, we prove \autoref{thm:inter-to-menger}; asserting that the infinite matroid intersection conjecture implies the infinite Menger theorem.

Given a graph $G$ and $S,T \subseteq V(G)$, a set $X \subseteq V(G)$ is called an $S$--$T$ \emph{separator} if $G-X$ contains no $S$--$T$ path.   
The infinite Menger theorem reads as follows.

\begin{thm}[Aharoni and Berger \cite{AharoniBerger}]\label{thm:menger}
Let $G$ be a connected graph.
Then for any $S,T\subseteq V(G)$ there is a set $\L$ of vertex-disjoint $S$--$T$ paths and an $S$--$T$ separator $X \subseteq \bigcup_{P\in \L} V(P)$
satisfying $|X \cap V(P)| = 1$ for each $P \in \L$. 
\end{thm}

Infinite matroid union cannot be used in order to obtain the infinite Menger 
Theorem directly via \autoref{thm:intersection} and \autoref{thm:inter-to-menger}.
Indeed, in \cite[Proposition 1.1]{union} we construct a finitary matroid $M$ and a co-finitary matroid $N$ such that their union is not a matroid.
Consequently, one cannot apply \autoref{thm:intersection} to the finitary matroids $M$ and $N^*$ in order to obtain \autoref{thm:nash-williams-aharoni} for them.
However, it is easy to see that \autoref{thm:nash-williams-aharoni} is true for these particular $M$ and $N^*$.

Next, we prove \autoref{thm:inter-to-menger}. 

\begin{proof}[Proof of \autoref{thm:inter-to-menger}.]
Let $G$ be a connected graph and let $S,T \subseteq V(G)$ be as in \autoref{thm:menger}.
We may assume that $G[S]$ and $G[T]$ are both connected.
Indeed, an $S$--$T$ separator with $G[S]$ and $G[T]$ connected gives rise to an $S$--$T$ separator when these are not necessarily connected.
Abbreviate $E(S):=E(G[S])$ and $E(T):=E(G[T])$, let $M$ be the finite-cycle matroid $M_F(G)$, and 
put $M_S:=M/E(S)-E(T)$ and $M_T:=M/E(T)-E(S)$; all three matroids are clearly finitary.

Assuming that the infinite matroid intersection conjecture holds for $M_S$ and $M_T$, there exists a set $I\in \I(M_S)\cap \I(M_T)$ which admits a partition $I=J_S \cup J_T$ 
satisfying
\[\cl_{M_S}(J_S)\cup \cl_{M_T}(J_T)=E,\]
where $E = E(M_S) = E(M_T)$.
We regard $I$ as a subset of $E(G)$.

For the components of $G[I]$ we observe two useful properties. 
As $I$ is disjoint from $E(S)$ and $E(T)$, the edges of a cycle in a component of $G[I]$ form a circuit in both, $M_S$ and $M_T$, contradicting the independence of $I$ in either. Consequently, 
\begin{equation}\label{eqn:x1}
 \text{the components of $G[I]$ are trees.}
\end{equation}
Next, an $S$-path\footnote{A non-trivial path meeting $G[S]$ exactly in its end vertices.} or a $T$-path in a component of $G[I]$ gives rise to a circuit of $M_S$ or $M_T$ in $I$, respectively.
Hence,
\begin{equation}\label{eqn:x2}
 \text{$|V(C)\cap S| \leq 1$ and $|V(C)\cap T| \leq 1$ for each component $C$ of $G[I]$.}
\end{equation}

Let $\C$ denote the components of $G[I]$ meeting both of $S$ and $T$.
Then by \eqref{eqn:x1} and \eqref{eqn:x2} each member of $\C$ contains a unique $S$--$T$ path and we denote the set of all these paths by $\L$.
Clearly, the paths in $\L$ are vertex-disjoint.

In what follows, we find a set $X$ comprised of one vertex from each $P\in\L$ to serve as the required $S$--$T$ separator. 
To that end, we show that one may alter the partition $I = J_S \cup J_T$ to yield a partition 
\begin{equation}\label{eqn:newpart}
	I = K_S \cup K_T\text{ satisfying } cl_{M_S}(K_S)\cup cl_{M_T}(K_T)=E \text{ and (Y.1-4)},
\end{equation}  
where (Y.1-4) are as follows. 
\begin{enumerate}[(Y.1)]
\item \label{y1} Each component $C$ of $G[I]$ contains a vertex of $S\cup T$.
\item \label{y2} Each component $C$ of $G[I]$ meeting $S$ but not $T$ satisfies $E(C) \subseteq K_S$.
\item \label{y3} Each component $C$ of $G[I]$ meeting $T$ but not $S$ satisfies $E(C) \subseteq K_T$.
\item \label{y4} Each component $C$ of $G[I]$ meeting both, $S$ and $T$, contains at most one vertex which at the same time
	\begin{enumerate}[(a)]
		\item \label{y4.a} lies in $S$ or is incident with an edge of $K_S$, and
		\item \label{y4.b} lies in $T$ or is incident with an edge of $K_T$.
	\end{enumerate}
\end{enumerate}

Postponing the proof of \eqref{eqn:newpart}, we first show how to deduce the existence of the required $S$--$T$ separator from \eqref{eqn:newpart}.
Define a pair of sets of vertices $(V_S,V_T)$ of $V(G)$ by letting $V_S$ consist of those vertices contained in $S$ or incident with an edge of $K_S$ and defining $V_T$ in a similar manner.  
Then $V_S \cap V_T$ may serve as the required $S$--$T$ separator.
To see this, we verify below  that $(V_S,V_T)$ satisfies all of the terms (Z.1-4) stated next.   
\begin{enumerate}[(Z.1)]
 \item $V_S\cup V_T=V(G)$;   \label{z1}
\item for every edge $e$ of $G$ either $e\subseteq V_S$ or $e\subseteq V_T$;  \label{z2}
\item every vertex in $V_S\cap V_T$ lies on a path from $\L$;  and \label{z3}
\item every member of $\L$ meets $V_S\cap V_T$ at most once.     \label{z4}
\end{enumerate}

To see (Z.\ref{z1}), suppose $v$ is a vertex not in $S\cup T$. As $G$ is connected, such a vertex is incident with some edge $e \notin E(T) \cup E(S)$. 
The edge $e$ is spanned by $K_T$ or $K_S$; say $K_T$.
Thus, $K_T+ e$ contains a circle containing $e$ or $G[K_T+e]$ has a $T$-path containing $e$. In either case $v$ is incident with an edge of $K_T$ and thus in $V_T$, as desired. 

To see (Z.\ref{z2}), let $e\in cl_{M_T}(K_T)\setminus K_T$; so that $K_T+e$ has a circle containing $e$ or $G[K_T+e]$ has $T$-path containing $e$; in either case both end vertices of $e$ are in $V_T$, as desired. The treatment of the case
$e\in cl_{M_S}(K_S)$ is similar. 

To see (Z.\ref{z3}), let $v\in V_S\cap V_T$; such is in $S$ or is incident with an edge of $K_S$, and in $T$ or is incident with an edge in $K_T$. Let $C$ be the component of $G[I]$ containing $v$. 
By (Y.1-4), $C\in\C$, i.e.\ it meets both, $S$ and $T$ and therefore contains an $S$--$T$ path $P\in\L$.
Recall that every edge of $C$ is either in $K_S$ or $K_T$ and consider the last vertex $w$ of a maximal initial segment of $P$ in $C-K_T$.
Then $w$ satisfies (Y.\ref{y4.a}), as well as (Y.\ref{y4.b}), implying $v=w$;
so that $v$ lies on a path from $\L$. 

To see (Z.\ref{z4}), we restate (Y.\ref{y4}) in terms of $V_S$ and $V_T$:
each component of $\C$ contains at most one vertex of $V_S\cap V_T$.
This clearly also holds for the path from $\L$ which is contained in $C$.

It remains to prove \eqref{eqn:newpart}.  
To this end, we show that any component $C$ of $G[I]$ contains a vertex of $S\cup T$.
Suppose not.
Let $e$ be the first edge on a $V(C)$--$S$ path $Q$ which exists by the connectedness of $G$.
Then $e\notin I$ but without loss of generality we may assume that $e\in\cl_{M_S}(J_S)$.
So in $G[I] + e$ there must be a cycle or an $S$-path.
The latter implies that $C$ contains a vertex of $S$ and the former means that $Q$ was not internally disjoint to $V(C)$, yielding contradictions in both cases.

We define the sets $K_S$ and $K_T$ as follows.
Let $C$ be a component of $G[I]$.
\begin{enumerate}
	\item If $C$ meets $S$ but not $T$, then include its edges into $K_S$.
	\item If $C$ meets $T$ but not $S$, then include its edges into $K_T$.
	\item Otherwise ($C$ meets both of $S$ and $T$) there is a path $P$ from $\L$ in $C$.
		Denote by $v_C$ the last vertex of a maximal initial segment of $P$ in $C-J_T$.
		As $C$ is a tree, each component $C'$ of $C-v_C$ is a tree and there is a unique edge $e$ between $v_C$ and $C'$.
		For every such component $C'$, include the edges of $C'+e$ in $K_S$ if $e\in J_S$ and in $K_T$ otherwise, i.e.\ if $e\in J_T$.  
\end{enumerate}
Note that, by choice of $v_C$, either $v_C$ is the last vertex of $P$ or the next edge of $P$ belongs to $J_T$.
This ensures that $K_S$ and $K_T$ satisfy (Y.\ref{y4}).
Moreover, they form a partition of $I$ which satisfies (Y.\ref{y1}-\ref{y3}) by construction.
It remains to show that $\cl_{M_S}(K_S)\cup \cl_{M_T}(K_T) = E$.

As $K_S\cup K_T = I$, it suffices to show that any $e \in E\setminus I$ is spanned by $K_S$ in $M_S$ or by $K_T$ in $M_T$.
Suppose $e\in \cl_{M_S}(J_S)$, i.e.\ $J_S + e$ contains a circuit of $M_S$.
Hence, $G[J_S]$ either contains an $e$-path $R$ or two disjoint $e$--$S$ paths $R_1$ and $R_2$.
We show that $E(R)\subseteq K_S$ or $E(R)\subseteq K_T$ in the former case and $E(R_1) \cup E(R_2) \subseteq K_S$ in the latter.

The path $R$ is contained in some component $C$ of $G[I]$.
Suppose $C\in \C$ and $v_C$ is an inner vertex of $R$.
By assumption, the edges preceding and succeeding $v_C$ on $R$ are both in $J_S$ and hence the edges of both components of $C-v_C$ which are met by $R$ plus their edges to $v_C$ got included into $K_S$, showing $E(R)\subseteq K_S$.
Otherwise $C\notin\C$ or $C\in \C$ but $v_C$ is no inner vertex of $R$.
In both cases the whole set $E(R)$ got included into $K_S$ or $K_T$.

We apply the same argument to $R_1$ and $R_2$ except for one difference.
If $C\notin\C$ or $C\in \C$ but $v_C$ is no inner vertex of $R_i$, then $E(R_i)$ got included into $K_S$ as $R_i$ meets $S$.

Although the definitions of $K_S$ and $K_T$ are not symmetrical, a similar argument shows $e\in \cl_{M_S}(K_S)\cup \cl_{M_T}(K_T)$ if $e$ is spanned by $J_T$ in $M_T$.
\end{proof}

Note that the above proof requires only that \autoref{thm:nash-williams-aharoni} holds for finite-cycle matroids.

\section{From infinite matroid union to infinite matroid intersection}\label{sec:union-to-inter}

In this section, we prove \autoref{thm:intersection}. 

\begin{proof}[Proof of \autoref{thm:intersection}.] Our starting point is the well-known proof from finite matroid theory that matroid union implies a solution to the matroid intersection problem. With that said, let $B_1 \cup B^*_2 \in \B(M_1 \vee M_2^*)$ where $B_1 \in \B(M_1)$ and $B^*_2 \in \B(M^*_2)$, and let $B_2 = E \setminus B^*_2 \in \B(M_2)$. Then, put $I = B_1 \cap B_2$ and note that $I \in \I(M_1) \cap \I(M_2)$. 
We show that $I$ admits the required partition. 

For an element $x \notin B_i$, $i=1,2$, we write $C_i(x)$ to denote the fundamental circuit of $x$ into $B_i$ in $M_i$.
For an element $x \notin B^*_2$, we write  $C^*_2(x)$ to denote the fundamental circuit of $x$ into $B^*_2$ in $M_2^*$. 
Put $X = B_1 \cap B^*_2$, $Y = B_2 \setminus I$, and $Z = B^*_2 \setminus X$, see \autoref{fig:conjecture}. 

\begin{figure} [htpb]   
\begin{center}
   	  \includegraphics[width=10cm]{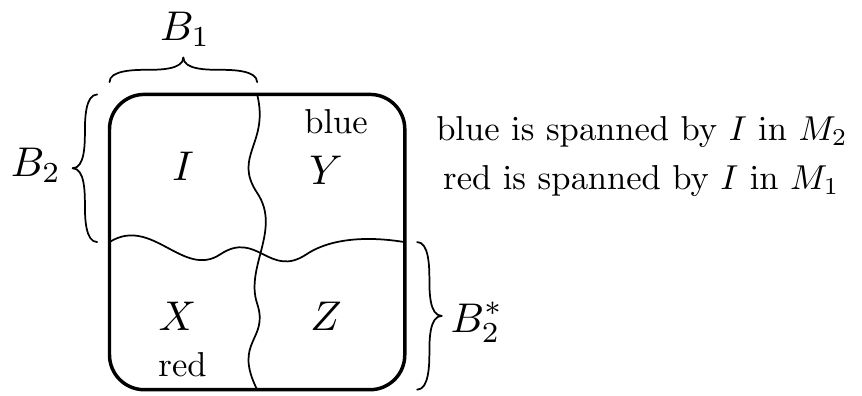}
   	  \caption{The sets $X$, $Y$, and $Z$ and their colorings.}
   	  \label{fig:conjecture}
\end{center}
\end{figure}

Observe that
\begin{equation}\label{eqn:span}
cl_{M_1}(I) \cup cl_{M_2}(I) = E = I\cup X \cup Y \cup Z.
\end{equation}
To see \eqref{eqn:span}, note first that
\begin{equation}\label{eqn:x}
X \subseteq cl_{M_2}(I).
\end{equation}
Clearly, no member of $X$ is spanned by $I$ in $M_1$. Assume then that $x \in X$ is not spanned by $I$ in $M_2$
 so that there exists a $y \in C_2(x) \cap Y$. Then, $x \in C^*_2(y)$, by \autoref{thm:meet}. Consequently, $B_1 \cup B_2^* \subsetneq B_1 \cup (B^*_2 + y -x) \in \I(M_1 \vee M_2^*)$; 
contradiction to the maximality of $B_1 \cup B^*_2$, implying \eqref{eqn:x}. 

By a similar argument, it holds that 
\begin{equation}\label{eqn:y}
Y \subseteq cl_{M_1}(I).
\end{equation}

To see that 
\begin{equation}\label{eqn:z}
Z \subseteq cl_{M_1}(I) \cup cl_{M_2}(I),
\end{equation}
assume, towards contradiction, that some $z \in Z$ is not spanned by $I$ neither in $M_1$ nor in $M_2$ so that 
there exist an $x \in C_1(z) \cap X$ and a $y \in C_2(z) \cap Y$. Then $B_1-x+z$ and $B_2-y+z$ are bases and thus 
$B_1 \cup B^*_2 \subsetneq (B_1-x+z) \cup (B^*_2-z+y)$; contradiction to the maximality of $B_1 \cup B^*_2$. Assertion \eqref{eqn:span} is proved. 

The problem of finding a suitable partition $I=J_1\cup J_2$ can be phrased as a (directed) graph coloring problem. 
By \eqref{eqn:span}, each $x\in E \setminus I$ satisfies $C_1(x) -x \subseteq I$ or $C_2(x) -x \subseteq I$.
Define $G=(V,E)$ to be the directed graph whose vertex set is $V = E \setminus I$ and whose edge set is given by 
\begin{equation}
E=\{(x,y): C_1(x) \cap C_2(y) \cap I \not= \emptyset\}.
\end{equation}
Recall that a \emph{source} is a vertex with no incoming edges and a \emph{sink} is a vertex with no outgoing edges.
As $C_1(x)$ does not exist for any $x \in X$ and $C_2(y)$ does not exist for any $y \in Y$,
it follows that 
\begin{equation}\label{eqn:source-sink}
\text{the members of $X$ are sinks and those of $Y$ are sources in $G$.}
\end{equation}
A $2$-coloring of the vertices of $G$, by say blue and red, is called \emph{divisive} if it satisfies the following:
\begin{enumerate}[(D.1)]
	\item $I$ spans all the blue elements in $M_1$;
	\item $I$ spans all the red elements in $M_2$; and
	\item $J_1\cap J_2=\emptyset$ where $J_1:=\left( \bigcup_{x\text{ blue}} C_1(x) \right) \cap I$ and $J_2:=\left( \bigcup_{x\text{ red}} C_2(x) \right) \cap I$.
\end{enumerate}
Clearly, if $G$ has a divisive coloring, then $I$ admits the required partition. 

We show then that $G$ admits a divisive coloring.
Color with blue all the sources.
These are the vertices that can only be spanned by $I$ in $M_1$.
Color with red all the sinks, that is, all the vertices that can only be spanned by $I$ in $M_2$.
This defines a partial coloring of $G$ in which all members of $X$ are red and those of $Y$ are blue.
Such a partial coloring can clearly be extended into a divisive coloring of $G$ provided that
\begin{equation}\label{eqn:forbidden}
\text{$G$ has no $(y,x)$-path with $y$ blue and $x$ red.}
\end{equation}
Indeed, given \eqref{eqn:forbidden} and \eqref{eqn:source-sink}, color all vertices reachable by a path from a blue vertex with the color blue, color all vertices from which a red vertex is reachable by a path with red, and color all remaining vertices with, say, blue. The resulting coloring is divisive. 

It remains to prove \eqref{eqn:forbidden}. We show that the existence of a path as in \eqref{eqn:forbidden} contradicts the following property:

\noindent
{\sl Suppose that $M$ and $N$ are matroids and $B \cup B'$ is maximal in $\I(M \vee N)$. Let $y\notin B \cup B'$ and let $x \in B \cap B'$. Then, (by~\autoref{thm:chain})
\begin{equation}\label{eqn:no-chain}
\text{there exists no $(B,B',y,x)$-chain;}
\end{equation}
}
(in fact, the contradiction in the proofs of \eqref{eqn:x},\eqref{eqn:y}, and \eqref{eqn:z} arose from simple instances of such forbidden chains).

Assume, towards contradiction, that $P$ is a $(y,x)$-path with $y$ blue and $x$ 
red; the intermediate vertices of such a path are not colored since they are not a sink nor a source. In what follows we use $P$ to construct a $(B_1,B^*_2,y_0,y_{2|P|})$-chain 
$(y_0,y_1,\ldots,y_{2|P|})$ such that $y_0 \in Y$, $y_{2|P|} \in X$, all odd indexed members of the chain are in $V(P) \cap Z$, 
and all even indexed elements of the chain other than $y_0$ and $y_{2|P|}$ are in $I$. Existence of such a chain would contradict \eqref{eqn:no-chain}.\vspace{\baselineskip}

\noindent
\textbf{Definition of $\boldsymbol{y_0}$.}
As $y$ is pre-colored blue then either $y \in Y$ or $C_2(y) \cap Y \not= \emptyset$. 
In the former case set $y_0 =y$ and in the latter choose $y_0 \in C_2(y) \cap Y$.\vspace{\baselineskip}

\noindent
\textbf{Definition of $\boldsymbol{y_{2|P|}}$.}
In a similar manner, $x$ is pre-colored red since either $x \in X$ or $C_1(x) \cap X \not= \emptyset$. In the former case, set $y_{2|P|} = x$ and in the latter case choose $y_{2|P|}\in C_1(x) \cap X$.\vspace{\baselineskip}

\noindent
\textbf{The remainder of the chain.}
Enumerate $V(P) \cap Z = \{y_1,y_3,\ldots,y_{2|P|-1}\}$ where the enumeration is with respect to the order of the vertices defined by $P$. Next, for an edge $(y_{2i-1},y_{2i+1}) \in E(P)$, let $y_{2i} \in C_1(y_{2i-1}) \cap C_2(y_{2i+1})\cap I$; such exists by the assumption that $(y_{2i-1},y_{2i+1}) \in E$. As $y_{2i+1} \in C^*_2(y_{2i})$ for all relevant $i$, by \autoref{thm:meet}, the sequence $(y_0,y_1,y_2,\ldots,y_{2|P|})$ is a $(B_1,B^*_2,y_0,y_{2|P|})$-chain
in $\I(M_1 \vee M^*_2)$.

This completes our proof of \autoref{thm:intersection}.
\end{proof}

Note that in the above proof, we do not use the assumption that $M_1 \vee M^*_2$ is a matroid; in fact,  we only need that $\I(M_1 \vee M^*_2)$ has a maximal element.

\section{The graphic nearly finitary matroids}\label{sec:graphic}

In this section we prove Propositions~\ref{thm:alg} and~\ref{thm:topo} yielding a characterization of the graphic nearly finitary matroids.

For a connected graph $G$, a maximal set of edges containing no finite cycles is called an \emph{ordinary spanning tree}. A maximal set of edges containing no finite cycles nor any double ray is called an \emph{algebraic spanning tree}. These are the bases of $M_F(G)$ and $M_A(G)$, respectively.  We postpone the discussion about $M_C(G)$ to \autoref{sec:topo}.

To prove Propositions~\ref{thm:alg} and~\ref{thm:topo}, we require the following theorem of Halin~\cite[Theorem 8.2.5]{DiestelBook10}.

\begin{thm}[Halin 1965]\label{thm:halin}
	If an infinite graph $G$ contains $k$ disjoint rays for every $k\in\N$, then $G$ contains infinitely many disjoint rays.
\end{thm}

\subsection{The nearly finitary algebraic-cycle matroids}\label{sec:alg-proof}

The purpose of this subsection is to prove \autoref{thm:alg}.
\begin{proof}[Proof of \autoref{thm:alg}.]
Suppose that $G$ has $k$ disjoint rays for every integer $k$;
so that $G$ has a set $\R$ of infinitely many disjoint rays by \autoref{thm:halin}.
We show that $M_A(G)$ is not nearly finitary.

The edge set of $\bigcup \R = \bigcup_{R \in \R} R$ is independent in $M_A(G)^{\fin}$ as it induces no finite cycle of $G$.
Therefore there is a base of $M_A(G)^{\fin}$ containing it;
such induces an ordinary spanning tree, say $T$, of $G$.
We show that
\begin{equation}\label{eqn:remove}
\text{$T - F$ contains a double ray for any finite edge set $F\subseteq E(T)$.}
\end{equation}
This implies that $E(T)\setminus I$ is infinite for every independent set $I$ of $M_A(G)$ and hence $M_A(G)$ is not nearly finitary.
To see \eqref{eqn:remove}, note that $T - F$ has $|F|+1$ components for any finite edge set $F\subseteq E(T)$ as $T$ is a tree and successively removing edges always splits one component into two.
So one of these components contains infinitely many disjoint rays from $\R$ and consequently a double ray.

Suppose next, that $G$ has at most $k$ disjoint rays for some integer $k$ and let $T$ be an ordinary spanning tree of $G$, that is, $E(T)$ is maximal in $M_A(G)^{\fin}$.
To prove that $M_A(G)$ is nearly finitary, we need to find a finite set $F\subseteq E(T)$ such that $E(T) \setminus F$ is independent in $M_A(G)$, i.e.\ it induces no double ray of $G$.
Let $\R$ be a maximal set of disjoint rays in $T$;
such exists by assumption and $|\R|\leq k$.
As $T$ is a tree and the rays of $\R$ are vertex-disjoint, it is easy to see that $T$ contains a set $F$ of $|\R| - 1$ edges such that $T-F$ has $|\R|$ components which each contain one ray of $\R$.
By maximality of $\R$ no component of $T-F$ contains two disjoint rays, or equivalently, a double ray.
\end{proof}

\subsection{The nearly finitary topological-cycle matroids}\label{sec:topo}

In this section we prove \autoref{thm:topo} that characterizes the nearly finitary topological-cycle matroids. Prior to that, we first define these matroids. 
To that end we shall require some additional notation and terminology on which more details can be found in~\cite{RD:HB:graphmatroids}.

An \emph{end of $G$} is an equivalence class of rays, where two rays are \emph{equivalent} if they cannot be separated by a finite edge set.
In particular, two rays meeting infinitely often are equivalent.
Let the \emph{degree} of an end $\omega$ be the size of a maximal set of vertex-disjoint rays belonging to $\omega$, which is well-defined~\cite{DiestelBook10}.
We say that a double ray \emph{belongs to} an end if the two rays which arise from the removal of one edge from the double ray belong to that end;
this does not depend on the choice of the edge.
Two rays that belong to the same end are called a \emph{topological cycle} of $G$.

For a graph $G$ the topological-cycle matroid of $G$, namely $M_C(G)$, has $E(G)$ as its ground set and its set of circuits consists of the finite and topological cycles.
In fact, every infinite circuit of $M_C(G)$ induces at least one double ray; provided that $G$ is locally finite~\cite{DiestelBook10}.

A graph $G$ has only finitely many disjoint rays if and only if $G$ has only finitely many ends, each with finite degree.
Also, note that 
\begin{equation}\label{eqn:deg2}
	\text{every end of a $2$-connected locally finite graph has degree at least $2$.}
\end{equation}
Indeed, applying Menger's theorem inductively, it is easy to construct in any $k$-connected graph for any end $\omega$ a set of $k$ disjoint rays of $\omega$.

Now we are in a position to start the proof of \autoref{thm:topo}.
\begin{proof}[Proof of \autoref{thm:topo}.]
If $G$ has only a finite number of vertex-disjoint rays then $M_A(G)$ is nearly finitary by \autoref{thm:alg}.
Since $M_A(G)^{\fin} = M_C(G)^{\fin}$ and $\I(M_A(G)) \subseteq \I(M_C(G))$, we can conclude that $M_C(G)$ is nearly finitary as well.

Now, suppose that $G$ contains $k$ vertex-disjoint rays for every $k \in \N$.
If $G$ has an end $\omega$ of infinite degree, then there is an infinite set $\R$  of vertex-disjoint rays belonging to $\omega$.
As any double ray containing two rays of $\R$ forms a circuit of $M_C(G)$, the argument from the proof of \autoref{thm:alg} shows that $M_C(G)$ is not nearly finitary.

Assume then that $G$ has no end of infinite degree.  
There are infinitely many disjoint rays, by \autoref{thm:halin}.
Hence, there is a countable set of ends $\Omega = \{\omega_1,\omega_2,\ldots\}$.

We inductively construct a set $\R$ of infinitely many vertex-disjoint double rays, one belonging to each end of $\Omega$.
Suppose that for any integer $n\geq 0$ we have constructed a set $\R_n$ of $n$ disjoint double rays, one belonging to each of the ends $\omega_1,\ldots,\omega_n$.
Different ends can be separated by finitely many vertices so there is a finite set $S$ of vertices such that $\bigcup \R_n$ has no vertex in the component $C$ of $G-S$ which contains $\omega_{n+1}$.
Since $\omega_{n+1}$ has degree~2 by \eqref{eqn:deg2}, there are two disjoint rays from $\omega_{n+1}$ in $C$ an thus also a double ray $D$ belonging to $\omega_{n+1}$.
Set $\R_{n+1}:= \R_n\cup \{D\}$ and $\R:=\bigcup_{n\in\N} \R_n$.

As $\bigcup \R$ contains no finite cycle of $G$, it can be extended to an ordinary spanning tree of $G$.
Removing finitely many edges from this tree clearly leaves an element of $\R$ intact.
Hence, the edge set of the resulting graph still contains a circuit of $M_C(G)$.
Thus, $M_C(G)$ is not nearly finitary in this case as well.
\end{proof}

In the following we shall propose a possible common generalization of Propositions \ref{thm:alg} and \ref{thm:topo} to all infinite matroids.
We call a matroid $M$ \emph{$k$-nearly finitary} if every base of its finitarization contains a base
of $M$ such that their difference has size at most $k$. 
Note that saying `at most $k$' is not equivalent to saying `equal to $k$', consider for example the algebraic-cycle matroid of the infinite ladder.
In terms of this new definition, Propositions \ref{thm:alg} and \ref{thm:topo} both state for a certain class of infinite matroids that each member of this class is $k$-nearly finitary for some $k$.
In fact, for all known examples of nearly finitary matroids, there is such a $k$.
This raises the following open question.

\begin{que}
Is every nearly finitary matroid $k$-nearly finitary for some $k$? 
\end{que}

\subsection{Graphic matroids and the intersection conjecture}

By
\autoref{thm:nearly-intersection},
 the intersection conjecture is true for $M_C(G)$ and $M_{FC}(H)$ for any two graphs $G$ and $H$ since the first is co-finitary and the second is finitary.
Using also \autoref{thm:topo}, we obtain the following.

\begin{cor}\label{thm:int_MT1}
Suppose that $G$ and $H$ are $2$-connected, locally finite graphs with only a finite number of vertex-disjoint rays.
Then, $M_C(G)$ and $M_C(H)$ satisfy the intersection conjecture.\qed
\end{cor}

Using \autoref{thm:alg} instead of \autoref{thm:topo}, we obtain the following.

\begin{cor}\label{thm:int_MT2}
Suppose that $G$ and $H$ are graphs with only a finite number of vertex-disjoint rays.
Then, $M_A(G)$ and $M_A(H)$ satisfy the intersection conjecture if both are matroids.\qed
\end{cor}

With a little more work, the same is also true for $M_{FC}(G)$, see \autoref{thm:int_MFC}.

\begin{proof}[Proof of \autoref{thm:int_MFC}.]
First we show that $(((M_C(G)^{\fin})^*)^{\fin})^*=M_C(G)$ if $G$ is locally
finite. 
Indeed, then $M_C(G)^{\fin}=M_{FC}(G)$,  $M_{FC}(G)^*$ is the matroids whose
circuits are the finite and infinite bonds of $G$, and its finitarization has as
its circuits the finite bonds of $G$.
And the dual of this matroid is $M_C(G)$, see \cite{matroid_axioms} for example.

Having showed that $(((M_C(G)^{\fin})^*)^{\fin})^*=M_C(G)$ if $G$ is locally
finite, we next show that if $M_C(G)$ is nearly finitary, then so is
$M_{FC}(G)^*$.
For this let $B$ be a base of $M_{FC}(G)^*$ and $B'$ be a base of
$(M_{FC}(G)^*)^{\fin}$.
Then $B'\setminus B=(E\setminus B)\setminus (E\setminus B')$. Now $E\setminus B$
is a base of $M_{FC}(G)=M_C(G)^{\fin}$ and
by the above  $E\setminus B'$ is a base of $M_C(G)$. Since $M_C(G)$ is nearly finitary, 
$B'\setminus B$ is finite, yielding that $M_{FC}(G)^*$ is nearly finitary. 

As $M_{FC}(G)^*$ is nearly finitary and $M_{FC}(H)$ is finitary, $M_{FC}(H)$ and $M_{FC}(G)$ satisfy the intersection conjecture by \autoref{thm:nearly-intersection}.
\end{proof}

A similar argument shows that if $G$ and $H$ are are $2$-connected, locally finite graphs with only a finite number of vertex-disjoint rays, then one can also prove that $M_{FC}(G)^*$ and $M_{FC}(H)^*$
satisfy the intersection conjecture. Similar results are true for $M_C(G)^*$ or $M_A(G)^*$ in place of $M_{FC}(G)^*$.

\bibliographystyle{plain}
\bibliography{literatur}

\end{document}